# The Berry-like Sentence in the First-order Peano Arithmetic System with the Operation of Factorial

T. Mei

(Central China Normal University, Wuhan, Hubei PRO, People's Republic of China
E-Mail: meitao@mail.ccnu.edu.cn    meitaowh@public.wh.hb.cn )

**Abstract:** A first-order Peano Arithmetical system with the operation of factorial (PAF) is introduced. For any formula $A(x)$ with a free variable $x$ in PAF, we define a corresponding *B-formula* $\exists x B_A(x)$ so as to there exists unique number that is smallest in all natural numbers satisfying the formula $A(x)$ that satisfies the formula $B_A(x)$ if $A(x)$ is satisfiable. And then, we construct a formula $\exists x B_D(x)$ which means that "there exists $x$, for any B-formula $\exists x B'_A(x)$ whose Gödel code is smaller than a constant $a$, $x$ does not satisfy $B'_A(x)$, and $x$ is the smallest in those numbers that have such character." However, $\exists x B_D(x)$ itself is a B-formula and its Gödel code is smaller than $a$. Thus, $\exists x B_D(x)$ is a version in PAF of the Berry sentence "The smallest positive integer not nameable in under eleven words" that itself is in only ten words .

## 1    Basic idea

A version of the Berry paradox[1, 2] arises from considering a sentence "The smallest positive integer not nameable in under eleven words." But the sentence itself is a specification for that number in only ten words.

The application of the Berry paradox in the information-theoretic approach had been discussed by G. J. Chaitin[3]. In this paper, we try to construct a Berry-like sentence in the first-order Peano Arithmetic system with the operation of factorial (PAF). The essential steps are as follows:

(1) For any formula $A(x)$ with a free variable $x$, we construct a corresponding formula $B_A(x)$ with a free variable $x$:

$$B_A(x):\ \ A(x) \wedge \forall y\, A(y) \rightarrow ((y \neq x) \rightarrow (y > x))$$

The above formula has many equivalent forms. However, we prescribe that the above standard form is only considered in spite of another equivalent forms.

Any formula with a free variable $x$ is false or satisfiable. And maybe many natural numbers satisfy a formula $A(x)$ if $A(x)$ is satisfiable; however, there exists unique number that is smallest in all natural numbers satisfying the formula $A(x)$ that satisfies the formula $B_A(x)$ if $A(x)$ is satisfiable. For example, any natural number satisfies the formula $x=x+0$; however, there is unique natural number (**0**) that satisfies the corresponding formula $(x=x+0) \wedge \forall y\, (y=y+0) \rightarrow ((y \neq x) \rightarrow (y > x))$.

For $B_A(x)$ which has the above standard form, we say that the *B-formula* $\exists x B_A(x)$ specifies a number $c$ that satisfies $B_A(c)$.

(2) Considering a 3-place relation $r(l, m, n)$ on the set of natural numbers, where $l$ means a Gödel code of a B-formula $\exists x B_A(x)$, $m$ means a Gödel code of a proof of the formula $\neg B_A(n)$, $B_A(n)$ is obtained by using the constant $n$ substitutes the free variable $x$ in $B_A(x)$. It can be proved that $r(l$, *m*, *n*) can be represented in PAF, we assume that the representation corresponding with *r*(*l*, *m*, *n*) in PAF is *R*(*x*, *y*, *z*).

(3) Considering a formula *D*(*x*) with a free variable *x*: $\forall y\ (y<a) \rightarrow \exists z\ R(y, z, x)$. *D*(*x*) means that, for any B-formula $\exists x B'_A(x)$ whose Gödel code *y* is smaller than a constant *a*, *x* does not satisfy $B'_A(x)$.

Note that those B-formulas whose Gödel code are smaller than *a* are finite; and, further, every B-formula $\exists x B_A(x)$ specifies just only one natural number if $B_A(x)$ is satisfiable, so the set of natural numbers specified by those B-formulas whose Gödel code are smaller than *a* are finite.

(4) For *D*(*x*), considering the corresponding $B_D(x)$: $D(x) \wedge \forall w\ D(w) \rightarrow ((w \neq x) \rightarrow (w > x))$. Note that the form of $B_D(x)$ accords with the standard form of $B_A(x)$ for any formula *A*(*x*) with a free variable *x*. And, further, the B-formula $\exists x B_D(x)$ means that "there exists *x*, for any B-formula $\exists x B'_A(x)$ whose Gödel code is smaller than a constant *a*, *x* does not satisfy $B'_A(x)$, and *x* is the smallest in those numbers that have such character." However, $\exists x B_D(x)$ itself is a B-formula; and, further, if we chose an appropriate number *a*, the Gödel code of the formula $\exists x B_D(x)$ can be smaller than *a*. Thus, a version in PAF of the Berry sentence "The smallest positive integer not nameable in under eleven words" is performed.

## 2 First-order Peano Arithmetical System with the Operation of Factorial (PAF)

In this section, based on the normal PA[4, 5], we state a standard structure of first-order Peano Arithmetic system by a slightly modified form. The main character of the modified form is that, apart from addition and multiplication, factorial is introduced as an essential operation. The advantage that factorial is introduced is that we can get bigger number by less symbols.

### 2.1 The symbols and the corresponding Gödel codes of PAF

All the symbols, corresponding names and Gödel codes of PAF are listed in Tab.1.

**Tab.1 The symbols of PAF**

| Name of symbol | | Symbol | Gödel code |
|---|---|---|---|
| bracket | left bracket | ( | 3 |
| | right bracket | ) | 5 |
| comma | | , | 7 |
| a constant symbol | | **0** | 9 |
| function symbols | successor | ' | 11 |
| | addition | + | 13 |
| | multiplication | × | 15 |
| | factorial | ! | 17 |
| predication symbols | equality | = | 19 |
| logical symbols | negation | ¬ | 21 |
| | implication | → | 23 |
| | universal quantifier | ∀ | 25 |
| variable symbols | | $x_0, x_1, x_2, ..., x_k, ...$ | $2k+27$, $k=0, 1, 2, 3, ...$ |

### 2.2 The terms of PAF

The terms of PAF are those finite sequences of symbols of PAF which satisfy the following

rules:

(T1) **0** is a term;

(T2) Every variable symbol $x_k$ is a term;

(T3) If $s$ is a term, then $s'$ is also a term;

(T4) If both $s$ and $t$ are terms, then the $s+t$ is also a term;

(T5) If both $s$ and $t$ are terms, then the $s \times t$ is also a term;

(T6) If $s$ is a term, then $s!$ is also a term;

(T7) Nothing else is a term.

## 2.3 The formulas and the corresponding Gödel codes of PAF

The formulas of PAF are those finite sequences of symbols of PAF which satisfy the following rules:

(F1) If both $s$ and $t$ are terms, then $s=t$ is a formula;

(F2) If $A$ is a formula, then $\neg A$ is also a formula;

(F3) If both $A$ and $B$ are formulas, then $A \to B$ is also a formula;

(F4) If $A$ is a formula and $x_k$ is one of variable symbols, then $\forall x_k(A)$ is also a formula;

(F5) Nothing else is a formula.

If a formula $F$ consists of a sequence of symbols $\alpha_0 \alpha_1 \alpha_2 ... \alpha_n$, and the Gödel code of $\alpha_k$ is $\#\alpha_k$, $k=1, 2, ..., n$, then the Gödel code $\#F$ of $F$ is $\#F = p_0^{\#\alpha_0} p_1^{\#\alpha_1} p_2^{\#\alpha_2} ... p_n^{\#\alpha_n}$, where $p_k$ is the $(k+1)$th prime number, $p_0 = 2, p_1 = 3, p_2 = 5, ...$.

## 2.4 The axiom schema of PAF

All the axiom schema and corresponding serial numbers of PAF are listed in Tab.2.

Tab.2 The axiom schema of PAF

| Serial numbers | Axiom schema | Explanatory note |
|---|---|---|
| 10 | $(A \to (B \to A))$ | Both $A$ and $B$ are any formulas |
| 11 | $(A \to (B \to C)) \to ((A \to B) \to (A \to C))$ | All $A$, $B$ and $C$ are any formulas |
| 12 | $(\neg A \to \neg B) \to (B \to A)$ | Both $A$ and $B$ are any formulas |
| 13 | $\forall x_k A \to A$ | $A$ is a formula, variable $x_k$ does not occur free in $A$ |
| 14 | $\forall x_k A[x_k] \to A[t]$ | $A[x_k]$ is a formula, variable $x_k$ occurs free in $A[x_k]$, $t$ is any term free for $x_k$ in $A[x_k]$ |
| 15 | $\forall x_k(A \to B) \to (A \to \forall x_k B)$ | $B$ is any formula; $A$ is a formula, variable $x_k$ does not occur free in $A$ |
| 16 | $x_k = x_k$ | $x_k$ is any variable |
| 17 | $(x_i = x_j) \to (x_j = x_i)$ | Both $x_i$ and $x_j$ are any variables |
| 18 | $(x_i = x_j) \to ((x_i = x_k) \to (x_j = x_k))$ | All $x_i$, $x_j$ and $x_k$ are any variables |

| 19 | $(x_i=x_j)\rightarrow(x_i'=x_j')$ | Both $x_i$ and $x_j$ are any variables |
|---|---|---|
| 20 | $\neg x_k'=0$ | $x_k$ is any variable |
| 21 | $(x_i'=x_j')\rightarrow(x_i=x_j)$ | Both $x_i$ and $x_j$ are any variables |
| 22 | $x_k+0=x_k$ | $x_k$ is any variable |
| 23 | $x_i+x_j'=(x_i+x_j)'$ | Both $x_i$ and $x_j$ are any variables |
| 24 | $x_k\times 0=0$ | $x_k$ is any variable |
| 25 | $x_i\times x_j'=x_i\times x_j+x_i$ | Both $x_i$ and $x_j$ are any variables |
| 26 | $0!=0'$ | |
| 27 | $(x_k')!=x_k'\times x_k!$ | $x_k$ is any variable |
| 28 | $(A[0]\wedge \forall x_k(A[x_k]\rightarrow A[x_k']))\rightarrow A[x_k]$ | $A[x_k]$ is a formula, variable $x_k$ occurs free in $A[x_k]$ |

Note that neither of the two symbols "[" and "]" are the symbols of PAF, $A[x_k]$ means only that the variable $x_k$ occurs in the formula $A$ of PAF.

**2.5 The deductive rules of PAF**

Both the deductive rules, corresponding serial numbers and the names of PAF are listed in Tab.3.

Tab.3  The deductive rules of PAF

| Serial numbers | Name | Deductive rule | Explanatory note |
|---|---|---|---|
| 29 | Modus Ponens (MP) | $\{A, A\rightarrow B\}\vdash B$ | Both $A$ and $B$ are any formulas |
| 30 | Generalization (Gen) | $A \vdash \forall x_k A$ | $A$ is any formula, $x_k$ is any variable |

**2.6 A proof of a formula and the corresponding Gödel codes in PAF**

A proof of a formula $F$ in PAF is those finite sequences of formulas $F_0, F_1, F_2, \ldots, F_n$ of PAF which satisfy the following conditions:

(1) For every formula $F_k$ ($k=0, 1, 2, \ldots, n$) in the sequences of $F_0, F_1, F_2, \ldots, F_n$,
    (i) $F_k$ is an axiom, or
    (ii) there are $i, j<k$ such that $F_k$ follows from $F_i$ and $F_j$ by the deductive rule MP, or
    (iii) there is $j<k$ such that $F_k$ follows from $F_j$ by the deductive rule Gen.

(2) $F_n$ is just the formula $F$ in the sequence of $F_0, F_1, F_2, \ldots, F_n$.

If the Gödel codes of the sequence of formulas $F_0, F_1, F_2, \ldots, F_n$ as a proof of a formula $F$ are $\#F_0, \#F_1, \#F_2, \ldots, \#F_n$, respectively, then the Gödel code $\#\{F_0, F_1, F_2, \ldots, F_n\}_{\text{PfF}}$ of the proof of the formula $F$ is: $\#\{F_0, F_1, F_2, \ldots, F_n\}_{\text{PfF}} = p_0^{\#F_0} p_1^{\#F_1} p_2^{\#F_2} \ldots p_n^{\#F_n}$, where $p_k$ is the $(k+1)$th prime number.

We see, PAF is composed of the normal PA and factorial as an original symbol and two axioms about the operation of factorial whose serial numbers in the Tab. 2 are 26 and 27, so PA is a subsystem of PAF.

## 3  The Construction of Berry-like Sentence in PAF

In this section, we use the standard symbols listed the Tab.1 to express formulas in PAF. For example, using "$\neg(F\rightarrow(\neg G))$" instead of "$F \wedge G$", "$\neg \forall x \neg$" instead of "$\exists x$", "$\forall u\neg(y=x+u+0')\rightarrow\neg \forall v\neg(x=y+v+0')$" instead of "$y \neq x$", "$\neg \forall w\neg(y=x+w+0')$" instead of "$y>x$"; And, in an actual formula, $x, y, u, v, w$, is one of $x_k$ ($k=0, 1, 2, 3, \ldots$), respectively.

For the 3-place relation $r(l, m, n)$ defined in the 1st section (2), we have

**Lemma 1** $r(l, m, n)$ is a recursive relation.

*Proof*: At first, we state an arithmetic to determine whether $r(l, m, n)$ is true:

① For a given natural number $l$, check whether $l$ is a Gödel code of a formula $F$ in PAF. If that are true, then

② Check whether $F$ has the standard B-formula form

$\neg \forall x(A(x)\rightarrow\neg \forall y A(y)\rightarrow(((\forall u\neg(y=x+u+0')\rightarrow\neg \forall v\neg(x=y+v+0'))\rightarrow(\neg \forall w\neg(y=x+w+0')))))$

where $x, y, u, v, w$, is one of $x_k$ ($k=0, 1, 2, 3, \ldots$), respectively.

Note that the above formula has many equivalent forms. However, we prescribe that the checkage is restricted whether the formula $F$ has the above standard form in spite of another equivalent forms. This restriction may be key for the feasibility of the arithmetic and does not affect the final result. Because we ask that the Berry-like sentence $\exists xB_D(x)$ that we shall construct has also such standard form.

If $F$ has the above standard B-formula form, then

③ For a given natural number $m$, check whether $m$ is a Gödel code of a prove of the formula

$A(0^{(n)})\rightarrow\neg \forall y A(y)\rightarrow(((\forall u\neg(y=0^{(n)}+u+0')\rightarrow\neg \forall v\neg(0^{(n)}=y+v+0'))\rightarrow(\neg \forall w\neg(y=0^{(n)}+w+0'))))$

where $0^{(n)}$ means $0''\ldots'$ (There are $n$ successors "$'$").

If that are true, then $r(l, m, n)$ is true for the given $l, m, n$, or false.

Now that an arithmetic to determine whether $r(l, m, n)$ is true is given, according to the Church's Thesis, $r(l, m, n)$ is recursive relation. The Lemma is proved.

Any recursive relation can be represented in PA[5]; And, further, PA is a subsystem of PAF, thus, any recursive relation can be represented in PAF. We therefore have

**Theorem 1** There exists a formula $R(x_0, x_1, x_2)$ to be as the representation corresponding with $r(l, m, n)$ in PAF so that

(i) If $r(l, m, n)$ is true, then PAF $\models R(0^{(l)}, 0^{(m)}, 0^{(n)})$;

(ii) If $r(l, m, n)$ is false, then PAF $\models \neg R(0^{(l)}, 0^{(m)}, 0^{(n)})$.

Assuming that the formula $R(x_{k-2}, x_{k-1}, x_k)$ is composed of $L_1$ standard symbols listed the Tab.1; And assuming that the maximal Gödel code in all $L_1$ symbols is $2k+27$ (which is corresponding with the symbol $x_k$), the formula $D(x_k)$ defined in the 1st section (3) is:

$\forall x_{k-2}((\neg \forall x_{k+1}\neg((\ldots((0''')!)!\ldots)!=x_{k-2}+x_{k+1}+0'))\rightarrow\neg \forall x_{k-1}\neg(R(x_{k-2}, x_{k-1}, x_k)))$.

Assuming that there are $L_2$ symbols of factorial, that means $a = (\ldots(3!)!\ldots)!$ (There are $L_2$ factorial

"!") in $D(x)$ defined in the 1st section (3), then $D(x_k)$ is composed of $L_1+3L_2+30$ standard symbols listed the Tab.1, the maximal Gödel code in all $L_1+3L_2+30$ symbols is $2(k+1)+27$ that is corresponding with the symbol $x_{k+1}$.

And, further, the formula $\exists x_k B_D(x_k)$ defined in the 1st section (4) is:

$\exists x_k B_D(x_k)$:  $\neg \forall x_k((D(x_k)) \to \neg \forall x_{k+2}((D(x_{k+2})) \to ((\forall x_{2k+4} \neg (x_{k+2}=x_k+x_{2k+4}+0') \to$

$\neg \forall x_{2k+5} \neg (x_k=x_{k+2}+x_{2k+5}+0')) \to (\neg \forall x_{2k+6} \neg (x_{k+2}=x_k+x_{2k+6}+0')))))$

Note that there are $k+2$ symbols $x_0, x_1, \ldots, x_k, x_{k+1}$ in the formula $D(x_k)$ and $k+2$ symbols $x_{k+2}, x_{k+3}, \ldots, x_{2k+2}, x_{2k+3}$ in the formula $D(x_{k+2})$. We see, the formula $\exists x_k B_D(x_k)$ is composed of $L=2(L_1+3L_2+30)+65=2L_1+6L_2+125$ standard symbols listed the Tab.1, and the maximal Gödel code in all $L$ symbols is $2(2k+6)+27$ (which is corresponding with the symbol $x_{2k+6}$).

Of course, in the actual formula $\exists x_k B_D(x_k)$, maybe we have to add some brackets, so we have $L=2L_1+6L_2+c$, where $c$ is a constant. Thus, the Gödel code of the formula $\exists x_k B_D(x_k)$ is

$$G=p_0^{\#\alpha_0} p_1^{\#\alpha_1} p_2^{\#\alpha_2} \ldots p_L^{\#\alpha_L} = 2^{\#\alpha_0} 3^{\#\alpha_1} 5^{\#\alpha_2} \ldots p_L^{\#\alpha_L},$$ and $\#\alpha_i \leq 2(2k+6)+27=4k+39$, $i=0,1,2,\ldots,L$.

Now we prove that $G<a=(\ldots(3!)!\ldots)!$ (There are $L_2$ factorial "!") by choosing an appropriate number $L_2$. At first, we have

$$G=p_0^{\#\alpha_0} p_1^{\#\alpha_1} p_2^{\#\alpha_2} \ldots p_L^{\#\alpha_L} = 2^{\#\alpha_0} 3^{\#\alpha_1} 5^{\#\alpha_2} \ldots p_L^{\#\alpha_L} < (2 \times 3 \times \ldots \times p_L)^{4k+39} < p_L^{(4k+39)(L+1)}.$$

**Lemma 2**[6] The $(n+1)$th prime number $p_n < 2^{2^{n+1}}$.

Choosing $L_2>2L_1+c+1$, then $L+1<7L_2$; and, further, choosing $L_2$ satisfies $2^{L_2}>7L_2(4k+39)$, from the Lemma 2 we have

$$G< p_L^{(4k+39)(L+1)} < (2^{2^{L+1}})^{(4k+39)(L+1)} = 2^{(4k+39)(L+1) \cdot 2^{L+1}} < 2^{7L_2(4k+39) \cdot 2^{7L_2}} < 2^{2^{8L_2}} < e^{e^{8L_2}},$$

where $e=2.71828\ldots$ is base of the natural logarithm.

**Lemma 3** (The Stirling Formula) $z! = \sqrt{2\pi z}\left(\dfrac{z}{e}\right)^z e^{\frac{\theta}{12z}}$, where $0<\theta<1$.

**Lemma 4** Assuming $a=(\ldots(3!)!\ldots)!$ (There are $N$ factorial "!"), then $a>e^{e^{8N}}$ when $N \geq 4$.

*Proof*: Set $a=b!$, $b=(\ldots(3!)!\ldots)!$ (There are $N-1$ factorial "!"), according to the Lemma 3,

$\ln a = b\ln\dfrac{b}{e} + \dfrac{1}{2}\ln(2\pi b) + \dfrac{\vartheta}{12b} > b$, because $b=(\ldots(3!)!\ldots)! \geq ((3!)!)!=(6!)!=720!$.

And, further, set $b=c!$, $c=(\ldots(3!)!\ldots)!$ (There are $N-2$ factorial "!"), according to the Lemma 3, $\ln\ln a > \ln b = c\ln\dfrac{c}{e} + \dfrac{1}{2}\ln(2\pi c) + \dfrac{\vartheta}{12c} > c$, because $c=(\ldots(3!)!\ldots)! \geq (3!)!=6!=720$.

On the other hand, it can be proved easily that $c>8N$ when $N \geq 4$, then we have $\ln\ln a > 8N$ when $N \geq 4$. The Lemma is proved.

Summarizing the above discussion, we have

**Theorem 2** If $L_2$ satisfies $L_2>4$, $L_2>2L_1+c+1$, $2^{L_2}>7L_2(4k+39)$, then the Gödel code of the formula $\exists x_k B_D(x_k)$ is smaller than $a=(\ldots(3!)!\ldots)!$ (There are $L_2$ factorial "!").

According to the analysis in the 1st section (4), the formula $\exists x_k B_D(x_k)$ satisfying the conditions stated in the Theorem 2 is a version in PAF of the Berry sentence "The smallest positive integer not nameable in under eleven words" which itself is in only ten words. However, we guess that $\exists x_k B_D(x_k)$ is false, that means that such $c$ which satisfies $B_D(c)$ does not exist. Because, if there existed such $c$, then a paradox would occur in PAF. And, whether PAF can prove $\exists x_k B_D(x_k)$, or whether PAF can prove $\neg \exists x_k B_D(x_k)$? Whether $\exists x_k B_D(x_k)$ can be undecidable in PAF? These questions will be discussed in further.